\documentclass[12pt]{amsart}
\usepackage{amsfonts,amsmath,amsthm,amssymb}
\usepackage{latexsym}
\usepackage{paralist}
\newtheorem{theorem}{Theorem}[section]

\newtheorem{lemma}[theorem]{Lemma}

\newtheorem{proposition}[theorem]{Proposition}
\newtheorem{observation}[theorem]{Observation}

\theoremstyle{definition}

\theoremstyle{remark}

\theoremstyle{remark}

\newcommand{\beql}[1]{\begin{equation}\label{#1}}
\newcommand{\eeq}{\end{equation}}

\newcommand \ov {\overline}

\newcommand{\R}{\ensuremath{\mathbb{R}}}

\newcommand{\Co}{\ensuremath{\mathbb{C}}}

\begin{document}

\title{Towards a classification of $6\times 6$ complex Hadamard matrices}

\author{M\'at\'e Matolcsi, Ferenc Sz\"oll\H{o}si}

\date{January, 2007.}

\address{M\'at\'e Matolcsi: Alfr\'ed R\'enyi Institute of Mathematics,
Hungarian Academy of Sciences, POB 127, H-1364 Budapest,
Hungary.}\email{matomate@renyi.hu}

\address{Ferenc Sz\"oll\H{o}si: Budapest University of Technology and Economics (BUTE), H-1111, Egry J. u. 1, Budapest, Hungary.}\email{szoferi@math.bme.hu}

\thanks{M. Matolcsi was supported by OTKA-T047276, T049301, PF64061.}

\begin{abstract}
Complex Hadamard matrices have received considerable attention in
the past few years due to their appearance in quantum information
theory. While a complete characterization is currently available
only up to order 5 (in \cite{haagerup}), several new constructions
of higher order matrices have appeared recently \cite{dita, karol,
BN, MM, sz}. In particular, the classification of {\it
self-adjoint} complex Hadamard matrices of order $6$ was completed
by Beuachamp and Nicoara in \cite{BN}, providing a previously
unknown non-affine one-parameter orbit. In this paper we classify
all \emph{dephased, symmetric} complex Hadamard matrices
\emph{with real diagonal} of order $6$. Furthermore, relaxing the
condition on the diagonal entries we obtain a new non-affine
one-parameter orbit connecting the Fourier matrix $F_6$ and
Di\c{t}\u{a}'s matrix $D_6$. This answers a recent question of
Bengtsson \& al. in \cite{BB}.
\end{abstract}

\maketitle

{\bf 2000 Mathematics Subject Classification.} Primary 05B20,
secondary 46L10.

{\bf Keywords and phrases.} {\it Complex Hadamard matrices}

\section{Introduction}\label{sec:intro}
Throughout this paper we will use the notation of \cite{karol,
web} for well-known complex Hadamard matrices such as $F_6$,
$D_6$, $B_6$ etc.

Original interest in {\it complex} Hadamard matrices arose in
connection with orthogonal pairs of maximal Abelian
$\ast$-subalgebras (MASA's) of the $n\times n$ matrices
\cite{haagerup, popa, HJ, MW, nic1}. Subsequently, it was realized
in \cite{wer} that complex Hadamard matrices also play an
essential role in constructions of teleportation and dense coding
schemes in quantum information theory. This fact has given a new
boost to the study of complex Hadamard matrices in recent years.
On the one hand, several new and general constructions of such
matrices have appeared \cite{dita, karol, MM, sz}. On the other
hand it is natural to try to {\it fully classify} complex Hadamard
matrices of small order, as such characterization is currently
available only up to order 5 in \cite{haagerup}. Recently some
progress has been made in the $6\times 6$ case in \cite{BN} where
all {\it self-adjoint} complex Hadamard matrices are
characterized, and in \cite{BB} where numerical evidence is given
of the existence of a conjectured 4-parameter family. While an
algebraic form of such a 4-parameter family (if it exists at all)
remains out of reach, in this paper we present a {\it previously
unknown non-affine one-parameter family} of $6\times 6$ complex
Hadamard matrices which {\it connects the Fourier matrix $F_6$ and
Di\c{t}\u{a}'s matrix $D_6$}. This result complements the recent
catalogue \cite{karol} and answers a question of \cite{BB},
proving that apart from the isolated matrix $S_6$ the set of known
$6\times 6$ Hadamard matrices is {\it connected}.

It is also important to mention that the $6\times 6$ case is
distinguished as the smallest dimension where the maximum number
of mutually unbiased bases (MUBs) is not known. It is well-known
that if $d$ is a prime power than the maximal number of MUBs in
$\Co^d$ is $d+1$. The existence of MUBs is equivalent to the
existence of $d\times d$ complex Hadamard matrices satisfying
certain conditions (see e.g. \cite{karol}). For the current status
of MUB-related problems and, in particular, the case $d=6$  we
refer to \cite{BB} and references therein. The recent discovery of
the new family in \cite{BN}, and the results of this paper may
well be useful in the resolution of the MUB problem in dimension
6, and may give an indication to the maximal number of MUB's in
dimension $d=pq$.

Throughout the paper we restrict attention to {\it dephased,
symmetric} complex Hadamard matrices of order 6 (the standard
terminology dephased meaning the normalization condition that all
entries of the first row and column are +1). It is quite natural
to study the symmetric case for two reasons. First, the results of
\cite{BN} show that it is hopeful to obtain closed algebraic
expressions if we require the matrix to satisfy certain symmetry
assumptions (in \cite{BN} the self-adjoint case was classified).
Second, the inspection of known $6\times 6$ complex Hadamard
matrices shows that many of them, such as $F_6$, $D_6$, $C_6$,
$S_6$, are equivalent to a symmetric one (throughout the paper we
use the standard notion of equivalence (see e.g. \cite{karol}),
i.e. $H_1$ and $H_2$ are equivalent, $H_1\cong H_2$, if
$H_1=D_1P_1H_2P_2D_2$ with unitary diagonal matrices $D_1, D_2$
and permutation matrices $P_1, P_2$). These two facts suggested
that the set of symmetric Hadamard matrices of order 6 is on the
one hand 'small' enough to be described in algebraic form and, on
the other hand, 'rich' enough to contain interesting families of
matrices. However, this intuitive approach turned out to be a
little too optimistic in the first respect, and we needed to put
further restrictions on the diagonal elements so that our
algebraic calculations come to a comprehensible end. Accordingly,
the outline of the paper is as follows. In Section \ref{sec2} we
fully classify {\it dephased, symmetric} complex Hadamard matrices
of order 6 {\it with real diagonal}. It turns out that under this
restriction well-known matrices emerge only. Therefore, in Section
\ref{sec3} we relax the condition on some diagonal entries and
this leads to the discovery of a new non-affine one-parameter
family.

\section{Symmetric matrices with real diagonal}\label{sec2}

First we recall a simple but extremely useful result of
\cite{haagerup} (see also \cite[Lemma 2.6]{BN}).

\begin{lemma}\label{trukk}
Let $u,v,s,t$ be complex numbers on the unit circle.
 Then \\
$(u+v)(\overline{s}+\overline{t})(\ov us+\ov vt)\in \R .$ $\hfill
\square$
\end{lemma}

We will also need the following elementary facts. In a dephased
Hadamard matrix the sum of the entries in each row is 0 (except
for the first row where the sum is, of course, $n$). Given a row
$(x_1,x_2, x_3, x_4, x_5, x_6)$ we distinguish two possibilities.
First, if $\Sigma :=\frac{x_1+x_2+x_3+x_4}{2}=0$ then $x_5=-x_6$.
Second, and more importantly, if $\Sigma\ne 0$ and $|\Sigma |\le
1$ then the last two coordinates are determined (up to change of
order) as
\begin{equation}\label{kov}x_{5,6}=-\Sigma\pm \textbf{i}\frac{\Sigma}{|\Sigma
|}\sqrt{1-|\Sigma |^2}
\end{equation}
The point is that $-2\Sigma=x_5+x_6$ and it is easy to see
geometrically that $x_5$ and $x_6$, being unit vectors, are
determined as above.

 The main
result of this section is the following

\begin{theorem}\label{thm1}
Let $H$ be a dephased, symmetric complex Hadamard matrix of order
$6$ with real diagonal. Then $H$ is equivalent to $S_6$ or $D_6$.
\end{theorem}

The proof is based on Lemma \ref{trukk}, and some considerations
similar to those in \cite{BN}.

The diagonal elements of $H$ belong to $\{-1, 1\}$ by assumption.
It is clear that there are either at least four 1's in the
diagonal or at most three. Therefore, after a possible permutation
of the rows and columns it is enough to consider the following two
possibilities for the diagonal of $H$:

\begin{equation}\label{diag}Diag(H)\in \{(1, 1, 1, 1, \ast,\ast), (1,-1, -1, -1,
\ast,\ast)\},\end{equation} where the $\ast$'s stand for $\pm 1$.

\begin{lemma}\label{lems6}
Let $H$ be a $6\times 6$ symmetric complex Hadamard matrix of the
form:
\begin{equation}H=\left[\begin{array}{cccccc}
1 & 1 & 1 & 1 & 1 & 1\\
1 & 1 & x & y & \ast & \ast\\
1 & x & 1 & z & \ast & \ast\\
1 & y & z & 1 & \ast & \ast\\
1 & \ast & \ast & \ast & \ast & \ast\\
1 & \ast & \ast & \ast & \ast & \ast\\
\end{array}\right]\end{equation}
Then
\begin{itemize}
\item[(a)] two of $x,y,z$ must be equal.
\item[(b)] $H$ is equivalent to $S_6$.
\end{itemize}
(Note that we do not assume here that the last two diagonal
entries are real; it is already implied by the above form.)
\end{lemma}

\begin{proof}
First we prove (a). Let us  denote
$h_{2,5}=u,h_{2,6}=v,h_{3,5}=s,h_{3,6}=t$. We will use Haagerup's
idea as in Lemma \ref{trukk}. By the orthogonality relations of
rows 1, 2, 3 we have
\begin{equation}2 + x + y = -(u + v)\end{equation}
\begin{equation}2 + \overline{x} + \overline{z} = -(\overline{s} + \overline{t})\end{equation}
\begin{equation}1 + x + \overline{x} + z \overline{y} = -(s \overline{u} + t \overline{v}).\end{equation}

Now, Lemma \ref{trukk} implies

\begin{equation}\label{1}
(2 + x + y)(2 + \overline{x} + \overline{z})(1 + x +
\overline{x} + z\overline{y}) \in \mathbb{R}. \end{equation}

By similar arguments we obtain

\begin{equation}\label{2}(2 +\overline{x} + \overline{y})(2 + y + z)(1 + \overline{y} + y
+ \overline{z} x)\in \mathbb{R}\end{equation}
\begin{equation}\label{3}(2 +
x + z)(2 + \overline{y} + \overline{z})(1 + z + \overline{z} + y
\overline{x})\in \mathbb{R}.\end{equation}

After summing up these three expressions and eliminating  real
terms we get

\begin{equation}\label{4}
\overline{x}^2y+\overline{x}y^2+x\overline{z}^2+x^2\overline{z}+\overline{y}^2z+\overline{y}z^2
+8\left(\overline{x}y+x\overline{z}+\overline{y}z\right)\in\mathbb{R}
\end{equation}

Since a complex number is real if and only if it is equal to its
conjugate, we obtain an equality if we replace each variable by
its conjugate (i.e. its reciprocal) in the above expression. The
resulting equality can then be rearranged by simple algebra to
yield

\begin{equation}\label{5}\left(x - y\right)\left(x - z\right)\left(y -
z\right)\left(xy+yz+zx + 8xyz + x^2yz + xy^2z +
xyz^2\right)=0\end{equation}

The last factor in the product is clearly non-zero by the triangle
inequality (one term has modulus 8, and the others have modulus
1). This proves $(a)$.

Now we turn to $(b)$. It is easy  to see that all $x=y$, $x=z$,
$y=z$ lead to equivalent Hadamard matrices by permutation, thus we
can assume without loss of generality that $x=y$.

Now substitute back to \eqref{1} and eliminate all real terms to
get

\begin{equation}\label{6}
x^2 + \overline{x}\overline{z}+\overline{x}z
\in\mathbb{R}.\end{equation} Therefore this expression equals its
conjugate and we obtain
\begin{equation}\label{7}x^2 + \overline{x}\overline{z}+\overline{x}z - \overline{x}^2-xz-x\overline{z}=0\end{equation}
which yields
\begin{equation}\label{8}\left(x^2-1\right)\left(z+x^2z-x-xz^2\right)=0\end{equation}

Here $x=y=1$ is clearly a contradiction, because the first two
rows of $H$ cannot be orthogonal.

To show that $x=y=-1$ is also impossible we need to consider two
subcases. If $x=y=-1=z$ then the rows of the leading $4\times 4$
minor of $H$ are already mutually orthogonal, therefore the last
two entries of the first four rows of $H$ should also be mutually
orthogonal, which is clearly impossible. If $x=y=-1\ne z$ then $H$
must be equivalent to a matrix of the following form

\beql{9} H=\left[\begin{array}{cccccc}
1 & 1 & 1 & 1 & 1 & 1\\
1 & 1 & -1 & -1 & u & -u\\
1 & -1 & 1 & z & -z & -1\\
1 & -1 & z & 1 & -1 & -z\\
1 & u & -z & -1 & \ast & \ast\\
1 & -u & -1 & -z & \ast & \ast\\
\end{array}\right].\eeq

The last two entries of rows 3, 4 are determined by \eqref{kov},
and the order of $-1$ and $z$ in the fourth row is forced by the
orthogonality of rows 3-4. The same orthogonality now implies
$z=\omega$ or $z=\omega^2$ (with $\omega$ being the third root of
unity). If $z=\omega$ then the orthogonality of rows 2 and 3
implies that $u=-\frac{\textbf{i}}{\sqrt3}$ which is a
contradiction. $z=\omega^2$ implies $u=\frac{\textbf{i}}{\sqrt3}$,
again a contradiction.

Therefore in \eqref{8} we must have
$\left(z+x^2z-x-xz^2\right)=(z-x)(1-xz)=0,$ which implies $z=x$ or
$z=\overline{x}$. The case $z=x=y$ leads again to contradiction
due to the following reasons. First, we cannot have $x=y=z=-1$ as
argued already above. Second, if $x=y=z\ne -1$ then \eqref{kov}
implies that two of the rows 2, 3, 4 must contain the same two
entries in the last two places. However, those two rows then
contain the same element in four places, therefore they cannot be
orthogonal.

Thus, we must have $y=x\neq \pm 1$, $z=\overline{x}$. In this case
$H$ must be equivalent with the following form
\beql{10}H=\left[\begin{array}{cccccc}
1 & 1 & 1 & 1 & 1 & 1\\
1 & 1 & x & x & \ast & \ast\\
1 & x & 1 & \overline{x} & u & v\\
1 & x & \overline{x} & 1 & v & u\\
1 & \ast & u & v & \ast & \ast\\
1 & \ast & v & u & \ast & \ast\\
\end{array}\right]\eeq
where $\Re[x]\le 0$ due to the orthogonality of rows 1 and 3.
After conjugating $H$ if necessary (and noting that $S_6$ is
equivalent to its conjugate) we can also assume that $\Im[x]\ge
0$. The third and fourth rows are being forced by \eqref{kov} as
follows (we are free to choose the order due to permutation
equivalence):
 \beql{u}u=-1-\Re[x]+
\textbf{i}\sqrt{-\Re[x]^2-2\Re[x]},\eeq
\beql{v}v=-1-\Re[x]-
\textbf{i}\sqrt{-\Re[x]^2-2\Re[x]}.\eeq
 Now,
due to the orthogonality of rows  3, 4 we get
$4\Re[x]^2+10\Re[x]+4=0$ which gives the only possible solution
$\Re[x]=-\frac{1}{2}$ and, by $\Im[x]\ge 0$, we get $x=\omega$.
Then by \eqref{u} and \eqref{v} $u=\omega ,$ $v=\omega^2$, and

\beql{s6elott}H=\left[\begin{array}{cccccc}
1 & 1 & 1 & 1 & 1 & 1\\
1 & 1 & \omega & \omega & \ast & \ast\\
1 & \omega & 1 & \omega^2 & \omega & \omega^2\\
1 & \omega & \omega^2 & 1 & \omega^2 & \omega\\
1 & \ast & \omega & \omega^2 & \ast & \ast\\
1 & \ast & \omega^2 & \omega & \ast & \ast\\
\end{array}\right]\eeq

All the remaining entries are determined uniquely (using
\eqref{kov} and orthogonality) and we obtain

\beql{s6}H=\left[\begin{array}{cccccc}
1 & 1 & 1 & 1 & 1 & 1\\
1 & 1 & \omega & \omega & \omega^2 & \omega^2\\
1 & \omega & 1 & \omega^2 & \omega & \omega^2\\
1 & \omega & \omega^2 & 1 & \omega^2 & \omega\\
1 & \omega^2 & \omega & \omega^2 & 1 & \omega\\
1 & \omega^2 & \omega^2 & \omega & \omega & 1\\
\end{array}\right]\eeq
This matrix is clearly equivalent to Tao's matrix $S_6$ (we note
here, that this matrix was published earlier in \cite{agaian},
page 104).
\end{proof}

Having classified dephased, symmetric Hadamard matrices with
diagonal $(1,1,1,1,\ast,\ast)$ we can  assume that the number of
1's in the diagonal are at most three, in other words, the
diagonal is $(1,-1,-1,-1,\ast,\ast)$.

\begin{lemma}\label{lemd6}
Let $H$ be a $6\times 6$ symmetric complex Hadamard matrix of the
form \beql{lemma2}H=\left[\begin{array}{cccccc}
1 & 1 & 1 & 1 & 1 & 1\\
1 & -1 & x & y & \ast & \ast\\
1 & x & -1 & z & \ast & \ast\\
1 & y & z & -1 & \ast & \ast\\
1 & \ast & \ast & \ast & \ast & \ast\\
1 & \ast & \ast & \ast & \ast & \ast\\
\end{array}\right]\eeq
Then $H$ is equivalent to $D_6$. (Note that we do not assume here
that the last two diagonal entries are real.)
\end{lemma}

\begin{proof}

We can assume that $\Im[x]\geq 0$ (by conjugating all entries of
$H$ if necessary, and noting that $D_6$ is equivalent to its
conjugate).

First assume that two of $x,y,z$ are equal, say $x=y$ (we are free
to choose due to permutation equivalence) and, furthermore,
$z=-x$. Then $H$ is equivalent to

\begin{equation}H=\left[\begin{array}{cccccc}
1 & 1 & 1 & 1 & 1 & 1\\
1 & -1 & x & x & -x & -x\\
1 & x & -1 & -x & u & -u\\
1 & x & -x & -1 & v & -v\\
1 & -x & u & v & \ast & \ast\\
1 & -x & -u & -v & \ast & \ast\\
\end{array}\right]\end{equation}

By the orthogonality of rows 3, 4

\begin{equation}\label{uesv} 1+1+\overline{x}+x+u\overline{v}+u\overline{v}=0, \end{equation}

which yields $u\overline{v}\in \mathbb{R}$, therefore $v=u$ or
$v=-u$. The former is not possible due to \eqref{uesv}, therefore
$v=-u$. Then, using again \eqref{uesv} we get $x+\overline{x}=0$
and hence, by the assumption $\Im[x]\ge 0$, $x=\textbf{i}$.

Then, using \eqref{kov} we obtain that the last two entries of
rows 5, 6 must be $-1$ and $\textbf{i}$ in some order. The same
order in both rows is not possible because then these rows would
agree in four entries and could not be orthogonal. Therefore we
have two possibilities:
\begin{equation}H=\left[\begin{array}{cccccc}
1 & 1 & 1 & 1 & 1 & 1\\
1 & -1 & \textbf{i} & \textbf{i} & -\textbf{i} & -\textbf{i}\\
1 & \textbf{i} & -1 & -\textbf{i} & u & -u\\
1 & \textbf{i} & -\textbf{i} & -1 & -u & u\\
1 & -\textbf{i} & u & -u & -1 & \textbf{i}\\
1 & -\textbf{i} & -u & u & \textbf{i} & -1\\
\end{array}\right], \ \ \ \ {\textrm{or}}\ \ \ \
H=\left[\begin{array}{cccccc}
1 & 1 & 1 & 1 & 1 & 1\\
1 & -1 & \textbf{i} & \textbf{i} & -\textbf{i} & -\textbf{i}\\
1 & \textbf{i} & -1 & -\textbf{i} & u & -u\\
1 & \textbf{i} & -\textbf{i} & -1 & -u & u\\
1 & -\textbf{i} & u & -u & \textbf{i} & -1\\
1 & -\textbf{i} & -u & u & -1 & \textbf{i}\\
\end{array}\right].\end{equation}

In the first case, by the orthogonality of rows $4, \ 5$ we get
$u=\pm\textbf{i}$. Both choices lead to Hadamard matrices
equivalent to $D_6$.

In the second case the orthogonality of rows 4, 5 yields $u=\pm
1$. Both choices lead to Hadamard matrices equivalent to $D_6$.

Let us now turn to the case when $x=y$, but $z\ne -x$. We will
show that it is not possible. The last two entries of row 3 must
be $-x$ and $-z$ and we are free to choose the order due to
permutation equivalence. Therefore,
\begin{equation}H=\left[\begin{array}{cccccc}
1 & 1 & 1 & 1 & 1 & 1\\
1 & -1 & x & x & -x & -x\\
1 & x & -1 & z & -x & -z\\
1 & x & z & -1 & \ast & \ast\\
1 & -x & -x & \ast & \ast & \ast\\
1 & -x & -z & \ast & \ast & \ast\\
\end{array}\right]\end{equation}

By the orthogonality of rows 2, 3 it  we get
\beql{11}1-\overline{x}-x+x\overline{z}+1+x\overline{z}=0,\eeq
thus $x\overline{z}\in \mathbb{R}$ implying $x=z$ or $x=-z$. The
first case is not possible by \eqref{11}, while the second
contradicts our current assumptions.

Now we turn to the case when $x,y,z$ are all distinct. Again, we
will show that it is not possible. We need to distinguish two
subcases. First we assume that one variable is the negative of
another, say $y=-x$ (we are free to choose due to permutation
equivalence). By assumption we cannot have $z=x$ or $z=-x$ (in the
latter case $z=y$ would hold). Therefore we have two choices to
fill up the last two entries of rows 3, 4:
\begin{equation}\label{1v2}H=\left[\begin{array}{cccccc}
1 & 1 & 1 & 1 & 1 & 1\\
1 & -1 & x & -x & u & -u\\
1 & x & -1 & z & -x & -z\\
1 & -x & z & -1 & x & -z\\
1 & u & -x & x & \ast & \ast\\
1 & -u & -z & -z & \ast & \ast\\
\end{array}\right] \ \ \ \textrm{or}\ \ \ \ H=\left[\begin{array}{cccccc}
1 & 1 & 1 & 1 & 1 & 1\\
1 & -1 & x & -x & u & -u\\
1 & x & -1 & z & -x & -z\\
1 & -x & z & -1 & -z & x\\
1 & u & -x & x & \ast & \ast\\
1 & -u & -z & -z & \ast & \ast\\
\end{array}\right]\end{equation}

In the first case the orthogonality of rows 3, 4 implies
$z=\pm\textbf{i}$ and we can assume (by conjugation if necessary)
that $z=\textbf{i}$. Then the orthogonality of rows $2, \ 3$ and
$2, \ 4$ yield the equalities
$1-x-\overline{x}-\overline{x}\textbf{i}-\overline{u}x+\overline{u}\textbf{i}=0$
and
$1+x+\overline{x}\textbf{i}+\overline{x}+\overline{u}x+\overline{u}\textbf{i}=0$
which, after summation, imply $u=-\textbf{i}$. Substituting back
$u=-\textbf{i}$ we get $x=\pm\textbf{i}$. But this implies $z=x$
or $z=y (=-x)$ which contradicts our current assumptions.

In the second case of \eqref{1v2} the orthogonality of rows 3, 4
implies $1-1-\overline{z}-z+x\overline{z}-\overline{x}z=0$, which
yields $x\overline{z}-\overline{x}z \in\mathbb{R}$. This is only
possible if $x\overline{z}\in\mathbb{R}$, that is $x=z$ or
$-x=z=y$, which are both excluded by assumption.

Lastly, assume that all $x,y,z$ are all distinct and none of them
is the negative of another. By formula \eqref{kov} there are the
following four different possibilities to fill out the last two
entries of rows 2, 3, 4 (we are free to fix the order of $-x$ and
$-y$ in row 2 due to permutation equivalence):

\begin{equation}\label{4eset}\left[\begin{array}{cc}
1 & 1\\
-x & -y\\
-x & -z\\
-y & -z\\
\ast & \ast\\
\ast & \ast\\
\end{array}\right],
\left[\begin{array}{cc}
1 & 1\\
-x & -y\\
-x & -z\\
-z & -y\\
\ast & \ast\\
\ast & \ast\\
\end{array}\right],
\left[\begin{array}{cc}
1 & 1\\
-x & -y\\
-z & -x\\
-y & -z\\
\ast & \ast\\
\ast & \ast\\
\end{array}\right],
\left[\begin{array}{cc}
1 & 1\\
-x & -y\\
-z & -x\\
-z & -y\\
\ast & \ast\\
\ast & \ast\\
\end{array}\right]\end{equation}
Now we analyze these cases separately.

CASE[1]: consider the first possibility listed in \eqref{4eset},
that is
\beql{kulxzy}H=\left[\begin{array}{cccccc}
1 & 1 & 1 & 1 & 1 & 1\\
1 & -1 & x & y & -x & -y\\
1 & x & -1 & z & -x & -z\\
1 & y & z & -1 & -y & -z\\
1 & -x & -x & -y & \ast & \ast\\
1 & -y & -z & -z & \ast & \ast\\
\end{array}\right]\eeq
Taking the inner product of rows 2, 3
\begin{equation}1-\overline{x}-x+y\overline{z}+1+y\overline{z}=0\end{equation}
Thus, $y\overline{z}\in\mathbb{R}$ and hence $z=y$ or $z=-y$ which
are both excluded by assumption.

CASE[2]: the second possibility listed in \eqref{4eset}. Consider
rows 2, 4 and apply the same simple argument as in CASE[1] above.

CASE[3]: the third possibility listed in \eqref{4eset}. The
orthogonality of rows 2, 3 implies
\begin{equation} 1-\overline{x}-x+y\overline{z}+x\overline{z}+\overline{x}y=0,\end{equation}
which means that $\overline{x}y+x\overline{z}+y\overline{z}\in
\mathbb{R}$. Therefore this expression equals its conjugate, i.e.
\begin{equation}\overline{x}y+x\overline{z}+y\overline{z}=x\overline{y}+\overline{x}z+\overline{y}z.\end{equation}
Using that conjugates are the same as reciprocals  this equation
is equivalent to
\begin{equation}(x+y)(x+z)(y-z)=0,\end{equation}
which contradicts our assumptions.

CASE[4]: the fourth possibility listed in \eqref{4eset}. Consider
rows 3, 4 and apply the same simple argument as in CASE[1] above.
\end{proof}

\section{A new family of $6\times 6$ complex Hadamard matrices}\label{sec3}

In Theorem \ref{thm1} we have classified all dephased symmetric
Hadamard matrices with real diagonal. Furthermore, from Lemmas
\ref{lems6} and \ref{lemd6} we see that it was, in fact, enough to
specify four real entries in the diagonal. It is then natural to
investigate the two remaining real options for the first four
entries of the diagonal, i.e. the cases $Diag H\in (1,-1,1,1,
\ast, \ast)$ and $Diag H\in (1,-1,-1,1,\ast,\ast)$.

Along the lines of Lemmas \ref{lems6} and \ref{lemd6} a
case-by-case argument shows that {\it there exists no dephased
symmetric complex Hadamard matrix with diagonal}
$(1,-1,-1,1,\ast,\ast)$. We do not include the details of this
fruitless calculation.

The last remaining case, $Diag H\in (1,-1,1,1, \ast, \ast)$, turns
out to be the most interesting one. Unfortunately we do not have a
full classification in this case, but we are able to obtain new
matrices nevertheless. For some preliminary calculations we
disregard the last +1 entry in the diagonal, and assume only that
$Diag H\in (1,-1,1,\ast, \ast, \ast)$. Due to the presence of the
$-1$ in the second row, the remaining entries must be $x,-x$ and
$y,-y$, and $H$ takes the form (up to permutation equivalence)

\begin{equation}H=\left[\begin{array}{rrrrrr}
1 & 1 & 1 & 1 & 1 & 1\\
1 & -1 & x & -y & y & -x\\
1 & x & 1 & a & b & c\\
1 & -y & a & \ast & \ast & \ast\\
1 & y & b & \ast & \ast & \ast\\
1 & -x & c & \ast & \ast & \ast\\
\end{array}\right].
\end{equation}

Using the orthogonality of rows 1, 2, 3 we apply Lemma \ref{trukk}
as follows.

\begin{equation}x+y=-\left(-y-x\right),\end{equation}
\begin{equation}2+\overline{x}+\overline{b}=-\left(\overline{a}+\overline{c}\right),\end{equation}
\begin{equation}1-x+\overline{x}+b\overline{y}=-\left(-a\overline{y}-c\overline{x}\right),\end{equation}
therefore
\begin{equation}\left(x+y\right)\left(2+\overline{x}+\overline{b}\right)\left(1-x+\overline{x}+b\overline{y}\right)\in\mathbb{R}
\end{equation}
After expanding and eliminating the real entries one gets:
\begin{equation}-2x^2-2xy+2\overline{x}y+\overline{x}^2y-x^2\overline{b}-xy\overline{b}+x\overline{y}b+y+b\in\mathbb{R}
\end{equation}
This expression therefore equals its conjugate and by simple
algebra we get
\begin{equation}
\left(y+x\right)\left[b^2(1+x^2)+b(2-x^3-2x^3y-2x^2-x+y+x^2y+2xy)-xy(1+x^2)\right]=0
\end{equation}
Unfortunately, we do not know how to handle the case when the
second factor equals zero, therefore we need to settle for the
{\it simplifying assumption} $y=-x$. We remark here, however, that
there exist non-trivial solutions with $y\ne -x$ too, such as the
permuted version of Bj\"orck's cyclic matrix
\begin{equation}C_6=\left[\begin{array}{cccccc}
 1 & 1 & 1 & 1 & 1 & 1 \\
 1 & -1 & d^2 & -d & d & -d^2 \\
 1 & d^2 & 1 & -d^3 & -\overline{d} & d^2 \\
 1 & -d & -d^3 & -d^3 & -d & d^4 \\
 1 & d & -\overline{d} & -d & \overline{d} & -1 \\
 1 & -d^2 & d^2 & d^4 & -1 & -d^4
\end{array}\right],\end{equation} where
$d=\frac{1-\sqrt{3}}{2}+\textbf{i}\cdot\sqrt{\frac{\sqrt{3}}{2}}.$

Having made the assumption $y=-x$ the matrix $H$ will now be
determined up to permutation equivalence and possible conjugation.
Now, $H$ takes the form
\begin{equation}H=\left[\begin{array}{rrrrrr}
1 & 1 & 1 & 1 & 1 & 1\\
1 & -1 & x & x & -x & -x\\
1 & x & 1 & a & b & c\\
1 & x & a & t & u & v\\
1 & -x & b & u & p & q\\
1 & -x & c & v & q & r\\
\end{array}\right]
\end{equation}

From the orthogonality of rows 2, 3 and 1, 3 we get:

\begin{equation}1-\overline{x}+x+x\overline{a}-x\overline{b}-x\overline{c}=0
\end{equation}
\begin{equation}\label{1-3}1+\overline{x}+1+\overline{a}+\overline{b}+\overline{c}=0
\end{equation}
Multiplying \eqref{1-3} by $x$ and then summing up and using
$\overline{x}=1/x$ we get
\begin{equation}a=\frac{x^2-2x-3}{2}.\end{equation}
This equation does have two solutions such that both $x$ and $a$
are on the unit ball,
\begin{equation}\label{x}x_{1,2}=\frac{1-\sqrt{13}}{3}\pm\textbf{i}\frac{\sqrt{-5+2\sqrt{13}}}{3}
\end{equation}
Let us take $x=x_1$ (the other choice leads to the conjugate
matrix), and hence
\begin{equation}\label{a}a=-\frac{7-\sqrt{13}}{9}-\textbf{i}\frac{\sqrt{19+14\sqrt{13}}}{9}\end{equation}
Now, since $2+x+a\neq0$ we can apply \eqref{kov} and obtain (up to
change of order, which we are free to choose due to permutation
equivalence)
\begin{equation}\label{b}b=\frac{-14+2\sqrt{13}-\sqrt{-58+34\sqrt{13}}}{18}-\textbf{i}\frac{\sqrt{134+22\sqrt{13}-8\sqrt{-2446+730\sqrt{13}}}}{18}
\end{equation}
\begin{equation}\label{c}c=\frac{-14+2\sqrt{13}+\sqrt{-58+34\sqrt{13}}}{18}+\textbf{i}\frac{\sqrt{134+22\sqrt{13}+8\sqrt{-2446+730\sqrt{13}}}}{18}
\end{equation}

Next we find $t,u$ and $v$. The orthogonality of rows 1, 4 and 2,
4 yield
\begin{equation}1+x+a+t+u+v=0,\end{equation}
\begin{equation}\label{ccc}1-x+a\overline{x}+t\overline{x}-u\overline{x}-v\overline{x}=0.\end{equation}
Multiplying \eqref{ccc} by $x$ and then summing up we obtain
\begin{equation}t=\frac{x^2-2x-1-2a}{2}.\end{equation}
Substituting the values of $x$ and $a$ we get $t=1$.

Then, using \eqref{kov} we obtain (the order being determined by
the orthogonality of rows 3, 4)
\begin{equation}u=c\end{equation}
\begin{equation}v=b\end{equation}

Finally we can use \eqref{kov} once again to complete rows 5 and 6
as (the order being determined by orthogonality of rows 4, 5, 6):
\begin{equation}\label{p}p=r=3-\sqrt{13}-\textbf{i}\sqrt{-21+6\sqrt{13}}\end{equation}
\begin{equation}\label{q}q=\frac{-19+4 \sqrt{13}}{9}+\textbf{i}\frac{2\sqrt{-122+38 \sqrt{13}}}{9}\end{equation}

Therefore, we have obtained
\begin{equation}M_6=\left[\begin{array}{rrrrrr}
1 & 1 & 1 & 1 & 1 & 1\\
1 & -1 & x & x & -x & -x\\
1 & x & 1 & a & b & c\\
1 & x & a & 1 & c & b\\
1 & -x & b & c & p & q\\
1 & -x & c & b & q & p\\
\end{array}\right],\end{equation}
where $x,a,b,c,p,q$ are determined by \eqref{x}, \eqref{a},
\eqref{b}, \eqref{c}, \eqref{p} and \eqref{q}, respectively. It is
easy to check that $M_6$ is indeed Hadamard. What we have shown
above is that {\it up to permutation equivalence $M_6$ and its
conjugate $M_6^\ast$ are the only dephased symmetric Hadamard
matrices with diagonal $(1,-1,1,\ast, \ast, \ast)$ and second row
consisting of the elements $(1,1,x,x,-x,-x)$}.

We will now proceed to show that $M_6$ is not contained in any of
the previously known $6\times 6$ families. We need the following
trivial
\begin{lemma}\label{selfadj} If a symmetric complex Hadamard matrix $H$
is equivalent to a self-adjoint one, then it is also equivalent to
its own conjugate i.e. $H\cong\overline{H}$.
\end{lemma}
\begin{proof}
Let $H$ be a symmetric complex Hadamard matrix and suppose that it
is equivalent to a self-adjoint one, say to $A=A^{\ast}$. Then
there are unitary diagonal $D_1,D_2$ and permutational matrices
$P_1, P_2$ such that
\begin{equation}P_1D_1HD_2P_2=A=A^{\ast}=P_2^{\ast}D_2^{\ast}H^{\ast}D_1^{\ast}P_1^{\ast}\end{equation}
By multiplying both sides with $D_2P_2$ from the left and $P_1D_1$
from the right we get:
\begin{equation}D_2P_2P_1D_1HD_2P_2P_1D_1=H^{\ast}=\overline{H}\end{equation}
This clearly says that $H \cong \overline{H}$.
\end{proof}

As a consequence we have
\begin{proposition}
$M_6$ and $M_6^\ast$ are not equivalent to any previously known
complex Hadamard matrix of order $6$.
\end{proposition}

\begin{proof}
We will use the Haagerup $\Lambda$-set of a matrix $H=[h_{jk}]$,
defined as
\begin{equation}\Lambda_H=\{ h_{ij}\overline h_{kj}h_{kl}\overline h_{il} \ \
{\textrm{for \ all}} \ \ 1\le i,j,k,l\le 6\}. \end{equation} It is
well-known that $\Lambda_H$ is invariant under equivalence (see
\cite{haagerup}).

The list of known $6\times 6$ Hadamard matrices is as follows. The
Fourier family $F_6^{(2)}(a,b)$ and its transposed $\left(
F_6^{(2)}(a,b)\right )^T$, the Di\c{t}\u{a} family $D_6^{(1)}(c)$,
and Tao's matrix $S_6$ are listed in \cite{karol}. The recently
discovered non-affine family $B_6^{(1)}(a)$ is given in \cite{BN}.

$M_6$ is clearly inequivalent to $S_6$ due to the Haagerup
$\Lambda$-set being different. $M_6$ is inequivalent to any matrix
in $F_6^{(2)}(a,b)$ since the third root of unity $\omega\in
\Lambda_{F_6^{(2)}(a,b)}$ for every matrix in that family (i.e.
for every $a,b$), while $\omega \notin \Lambda_{M_6}$. The same is
true for the transposed family $\left( F_6^{(2)}(a,b)\right )^T$.
$M_6$ is inequivalent to any matrix in $D_6^{(1)}(c)$ since
$\textbf{i}\in\Lambda_{D_6^{(1)}(c)}$ for every matrix in that
family, while $\textbf{i}\notin\Lambda_{M_6}$. Finally, using
Lemma \ref{selfadj}, $M_6$ is inequivalent to any of the matrices
contained in the self-adjoint non-affine family $B_6^{(1)}(a)$
since $M_6$ is inequivalent to $\overline{M}_6$, which can be seen
as follows. Let $a$ be as in $M_6$. One can easily check that
$\overline{a}^2\in \Lambda_{M_6}$ and $a^2\notin \Lambda_{M_6}$,
while $\overline{a}^2\notin \Lambda_{\overline{M}_6}$ and $a^2\in
\Lambda_{\overline{M}_6}$.

The same proof works for $M_6^\ast$, too.
\end{proof}

We  will now generalize our discrete result above and construct a
continuous one-parameter family stemming from $M_6$. As $M_6$ is
in block-matrix form it is quite natural to try and replace the
$1$'s on the main diagonal by some parameter $d$ and consider
dephased symmetric matrices of the following form:

\begin{equation}\label{fam}M_6(x)=\left[\begin{array}{rrrrrr}
1 & 1 & 1 & 1 & 1 & 1\\
1 & -1 & x & x & -x & -x\\
1 & x & d & a & b & c\\
1 & x & a & d & c & b\\
1 & -x & b & c & p & q\\
1 & -x & c & b & q & p\\
\end{array}\right],\end{equation}

The orthogonality of rows 2, 3 and 1, 3 imply

\begin{equation}\label{uj23}1-x+d\overline{x}+a\overline{x}-b\overline{x}-c\overline{x}=0\end{equation}
\begin{equation}1+x+d+a+b+c=0\end{equation}
Multiplying \eqref{uj23} by $x$ and then summing up we get
\begin{equation}\label{auj}a=\frac{x^2-2x-1}{2}-d.\end{equation}
It is easy to see that $0<|x^2-2x-1|<4$ for each choice of $x$ on
the unit circle, therefore we can use \eqref{auj} to apply
\eqref{kov} to obtain the values of $a$ and $d$ as follows (we are
free to choose the order due to permutation equivalence):
\begin{equation}\label{auj1}a=\frac{x^2-2x-1}{4}-\textbf{i}\frac{\left(x^2-2x-1\right)\sqrt{16-\left|x^2-2x-1\right|^2}}{4\left|x^2-2x-1\right|}
\end{equation}
\begin{equation}\label{duj}d=\frac{x^2-2x-1}{4}+\textbf{i}\frac{\left(x^2-2x-1\right)\sqrt{16-\left|x^2-2x-1\right|^2}}{4\left|x^2-2x-1\right|}\end{equation}

From \eqref{auj} we see that $1+x+a+d=\frac{x^2+1}{2}$ which
vanishes if and only if $x=\pm \textbf{i}$. We exclude $x=\pm
\textbf{i}$ from these considerations and remark that this case
can be handled separately as in Lemma \ref{lemd6}, and leads to
$H\cong D_6$. For $x\ne\pm \textbf{i}$ the values of $b$ and $c$
are determined uniquely  by \eqref{kov} as follows (we are free to
choose the order due to permutation equivalence):
\begin{equation}\label{buj}b=-\frac{1+x^2}{4}-\textbf{i}\frac{\left(1+x^2\right)\sqrt{16-\left|1+x^2\right|^2}}{4\left|1+x^2\right|}\end{equation}
\begin{equation}\label{cuj}c=-\frac{1+x^2}{4}+\textbf{i}\frac{\left(1+x^2\right)\sqrt{16-\left|1+x^2\right|^2}}{4\left|1+x^2\right|}\end{equation}

It is easy to check (rather by computer) that with these
parametric choices the first four rows of $H$ are mutually
orthogonal to each other.

We evaluate $1-x+b+c$ in order to use \eqref{kov} again to
determine $p$ and $q$. The orthogonality of rows 1, 5 and 2, 5
imply
\begin{equation}1-x+b+c+p+q=0\end{equation}
\begin{equation}\label{uj25}1+x+\overline{x}b+\overline{x}c-\overline{x}p-\overline{x}q=0.\end{equation}
Multiplying \eqref{uj25} by $x$ and summing up we get
\begin{equation}1-x+b+c=\frac{-x^2-2x+1}{2}\end{equation}
and we see that it does not vanish for unit vectors $x$. Therefore
we can apply \eqref{kov} to determine $p$ and $q$ as follows (the
order now being forced by orthogonality):
\begin{equation}\label{puj}p=\frac{ x^2+2x-1}{4}+\textbf{i}\frac{\left(x^2+2x-1\right)\sqrt{16-\left|x^2+2x-1\right|^2}}{4\left|x^2+2x-1\right|}\end{equation}
\begin{equation}\label{quj}q=\frac{x^2+2x-1}{4}-\textbf{i}\frac{\left(x^2+2x-1\right)\sqrt{16-\left|x^2+2x-1\right|^2}}{4\left|x^2+2x-1\right|}\end{equation}

It is easy to check by computer that the emerging matrix $H$ is
Hadamard for all $x\ne\pm \textbf{i}$. Therefore we have obtained
the following

\begin{theorem}
There is a one parameter \emph{symmetric} non-affine family of
Hadamard matrices given by \eqref{fam}, with $x=e^{\textbf{i}t},$
$x\ne\pm \textbf{i},$ and $a,b,c,d,p,q$ being given as in
\eqref{auj1}, \eqref{buj}, \eqref{cuj}, \eqref{duj}, \eqref{puj},
\eqref{quj}, respectively.
\end{theorem}

Finally we make the following observation which answers a question
raised in \cite{BB} and shows that the set of currently known
Hadamard matrices of order 6 is connected except for the isolated
matrix $S_6$. In particular, the family $M_6(x)$ above connects
the Fourier matrix $F_6$ and Di\c{t}\u{a}'s matrix $D_6$.

\begin{observation}
$M_6(1)\cong F_6$, $\lim\limits_{t\rightarrow
\frac{3\pi}{2}}M_6\left(e^{\textbf{i}t}\right)\cong D_6$, and
finally,
$M_6\left(e^{\textbf{i}\arccos\left(\frac{1-\sqrt{13}}{3}\right)}\right)=M_6$.
\end{observation}

\begin{equation}\lim\limits_{t\rightarrow
\frac{3\pi}{2}}M_6\left(e^{\textbf{i}t}\right)=\left[\begin{array}{rrrrrr}
 1 & 1 & 1 & 1 & 1 & 1 \\
 1 & -1 & -\textbf{i} & -\textbf{i} & \textbf{i} & \textbf{i} \\
 1 & -\textbf{i} & -1 & \textbf{i} & -1 & 1 \\
 1 & -\textbf{i} & \textbf{i} & -1 & 1 & -1 \\
 1 & \textbf{i} & -1 & 1 & -\textbf{i} & -1 \\
 1 & \textbf{i} & 1 & -1 & -1 & -\textbf{i}
\end{array}\right],\end{equation}

\begin{equation}M_6\left(1\right)=\left[\begin{array}{rrrrrr}
 1 & 1 & 1 & 1 & 1 & 1 \\
 1 & -1 & 1 & 1 & -1 & -1 \\
 1 & 1 & \omega^2 & \omega & \omega^2 & \omega \\
 1 & 1 & \omega & \omega^2 & \omega & \omega^2 \\
 1 & -1 & \omega^2 & \omega & -\omega^2 & -\omega \\
 1 & -1 & \omega & \omega^2 & -\omega & -\omega^2
\end{array}\right].\end{equation}

In summary, we have given a full classification of dephased
symmetric complex Hadamard matrices of order 6 with real diagonal,
showing that in this case the well-known matrices $S_6$ and $D_6$
emerge only. Furthermore, relaxing the reality condition on the
diagonal entries we have been able to obtain a new non-affine
family of dephased symmetric Hadamard matrices of order 6 which
connects $F_6$ and $D_6$.

It would be interesting to see whether this family can be extended
by further parameters. For example, it is natural to try to
replace the second entry of the diagonal ($-1$) by a parameter
$h$. However, currently we are unable to classify that case. It
also remains to be seen whether this new family helps to increase
the number of bases appearing in MUBs in $\Co^6$.

\end{document}